\documentclass[12pt]{amsart}
\usepackage{amscd,amsmath,amssymb,amsfonts}
\theoremstyle{plain}
\newtheorem{thm}{Theorem}
\newtheorem{lem}[thm]{Lemma}

\newtheorem{prop}[thm]{Proposition}

\theoremstyle{definition}
\newtheorem{defn}[thm]{Definition}
\newtheorem{rmk}[thm]{Remark}
\newtheorem{rmks}[thm]{Remarks}

\newtheorem{nota}[thm]{Notations}
\newtheorem{chal}[thm]{Challenge}
\numberwithin{thm}{section}
\numberwithin{equation}{section}

\newcommand{\eq}[2]{\begin{equation}\label{#1}#2 \end{equation}}
\newcommand{\ml}[2]{\begin{multline}\label{#1}#2 \end{multline}}
\newcommand{\ga}[2]{\begin{gather}\label{#1}#2 \end{gather}}

\newcommand{\surj}{\twoheadrightarrow}
\newcommand{\inj}{\hookrightarrow}

\newcommand{\Spec}{{\rm Spec \,}}

\newcommand{\sA}{{\mathcal A}}

\newcommand{\sI}{{\mathcal I}}

\newcommand{\sO}{{\mathcal O}}

\newcommand{\sR}{{\mathcal R}}
\newcommand{\sS}{{\mathcal S}}

\newcommand{\sZ}{{\mathcal Z}}
\newcommand{\A}{{\mathbb A}}

\newcommand{\G}{{\mathbb G}}

\renewcommand{\P}{{\mathbb P}}

\newcommand{\Z}{{\mathbb Z}}

\begin{document}

\title[Additive Chow groups]{An additive version of higher Chow
groups} 
\author{Spencer Bloch}
\address{Dept. of Mathematics,
University of Chicago,
Chicago, IL 60637,
USA}
\email{bloch@math.uchicago.edu}

\author{H\'el\`ene Esnault}
\address{Mathematik,
Universit\"at Essen, FB6, Mathematik, 45117 Essen, Germany}
\email{esnault@uni-essen.de}
\date{December 10, 2001}
\begin{abstract}  
The cosimplicial scheme 
$$\Delta^\bullet = \Delta^0 \begin{smallmatrix}
\to \\
\to \end{smallmatrix}   \Delta^1 \begin{smallmatrix}
\to \\ 
\to \\
\to \end{smallmatrix}
\ldots;\quad \Delta^n :=\Spec\Big(k[t_0,\dotsc,t_n]/(\sum t_i -t)\Big)$$ 
was used in \cite{B} to define higher Chow groups. In this note, we let $t$ 
tend to 0 and replace $\Delta^\bullet$  by a degenerate version 
$$Q^\bullet = Q^0 \begin{smallmatrix}
\to \\
\to \end{smallmatrix}  
 Q^1 
\begin{smallmatrix}
\to \\ 
\to \\
\to \end{smallmatrix}
 \ldots;\quad Q^n :=
\Spec\Big(k[t_0,\dotsc,t_n]/(\sum t_i)\Big)  $$ 
to define an additive version of the higher Chow groups. For a field $k$, we show
the Chow group of $0$-cycles on $Q^n$ in this theory is isomorphic to the absolute
$(n-1)$-K\"ahler forms $\Omega^{n-1}_k$. 

An analogous degeneration on the level of de Rham cohomology associated to
``constant modulus'' degenerations of varieties in various contexts is discussed. 
\end{abstract}
\subjclass{Primary 14C15, 14C35}
\maketitle
\begin{quote}

\end{quote}

\section{Introduction}

The purpose of this note is to study a common sort of limiting phenomenon which occurs
in the study of motives. Here is a simple example. Let $k$ be a field. Let
$S=\A^1_t=\Spec(k[t])$ and let
$T = \Spec(k[x,t]/(x(x-t)) \inj \A^2_{t,x}$. Over $S[1/t]$ the Picard scheme
$\text{Pic}(\A^2_{x,t},T)/S$ is represented by $\G_{m,S[1/t]}$. On the other hand, when
$t=0$ one gets $\text{Pic}(\A_x^1,\{x^2=0\})\cong \G_{a,k}$. In some sense, $\G_m$ has
``jumped'' to $\G_a$. 

In higher dimension, for $t\neq 0$, the Chow groups of the cosimplicial scheme
\eq{1.1a}{\Delta^\bullet = \Delta^0 \begin{smallmatrix}
\to \\ \to \end{smallmatrix}
 \Delta^1 \begin{smallmatrix} \to \\ \to \\ \to \end{smallmatrix}
\ldots;\quad \Delta^n :=
\Spec\Big(k[t_0,\dotsc,t_n]/(\sum t_i -t)\Big) 
}
are known to give motivic cohomology (\cite{V}). What can one say about the Chow groups of
the degenerate cosimplicial complex:
\eq{1.2a}{Q^\bullet = Q^0 \begin{smallmatrix} \to \\ \to
\end{smallmatrix}  Q^1  \begin{smallmatrix} \to \\ \to \\ \to
\end{smallmatrix}  \ldots;\quad Q^n :=
\Spec\Big(k[t_0,\dotsc,t_n]/(\sum t_i)\Big)  
}
Our main result is a calculation of the Chow groups of $0$-cycles on $Q^\bullet$. Let
$\sZ^n(Q^r)$ be the free abelian group on codimension $n$ algebraic cycles on $Q^r$
satisfying a suitable general position condition with respect to the face maps. Let
$SH^n(k,r)$ be the cohomology groups of the complex
\eq{1.3a}{\ldots \sZ^n(Q^{r+1}) \to \sZ^n(Q^{r}) \to \sZ^n(Q^{r-1}) \to \ldots
}
where the boundary maps are alternating sums of pullbacks along face maps.
Write $\Omega^\bullet_k$ for the absolute K\"ahler differentials.
\begin{thm}$SH^n(k,n) \cong \Omega^{n-1}_k$.
\end{thm}

Note for $n=1$, this is the above $\G_a$.

For another example of this limiting phenomenon, consider $\P^{n+1}$ with
coordinates
$U_0,\dotsc,U_{n+1}$ over a field $k$. Let
\eq{1.4a}{X:f(U_0,\dotsc,U_{n+1})=0
}
be homogeneous of degree
$n+2$. Let $V_i = t^{r_i}U_i,\ 0\le i\le n+1$, where $t$ is a
variable,
$r_0=0$, and $r_i \ge 0$. Let $N$ be the degree in $t$ of $F(t,U) :=
f(U_0,t^{r_1}U_1,\dotsc,t^{r_{n+1}}U_{n+1})$, and assume $s:=N-\sum r_i>0$. Write $u_i =
U_i/U_0,\ v_i = V_i/V_0$. Define a log $(n+1)$-form on
$\P^{n+1}\setminus X$
\eq{1.1}{\omega_f := \frac{du_1\wedge\ldots\wedge du_{n+1}}{f(1,u_1,\dotsc,u_{n+1})}.
}
We view $\omega_f$ (more precisely $pr_1^*\omega_f$) as a closed
$(n+1)$-form on
$(\P^{n+1} \setminus X)\times \G_{m,t}$. Making the substitution $u_i = t^{-r_i}v_i$ and clearing
denominators, the corresponding form {\it relative to } $\G_{m,t}$ can be written
\ml{1.2}{\omega_{f, {\rm rel}} = \frac{t^{-\sum r_i}dv_1\wedge\ldots\wedge
v_{n+1}}{f(1,t^{-r_1}v_1,\dotsc,t^{-r_{n+1}}v_{n+1})} =
\frac{t^{s}dv_1\wedge\ldots\wedge
v_{n+1}}{t^Nf(1,t^{-r_1}v_1,\dotsc,t^{-r_{n+1}}v_{n+1})} \\  =: 
t^{s}\nu_{f, {\rm rel}} }
Write $\nu_f$ for the evident absolute form lifting the relative
form $\nu_{f,{\rm rel}}$.
Define an $n$-form $\gamma_f$ by the equation
\eq{1.4}{\omega_f = t^s\nu_f + st^{s-1}dt\wedge\gamma_f
}
Then, assuming $s^{-1} \in k$
\eq{1.5}{\nu_f|_{t=0} = d\gamma_f|_{t=0}.
}

To combine cycles and forms, let
\eq{1.8}{ \Delta_{n+1}(1,u_1,\dotsc,u_{n+1}) := u_1\cdot u_2\cdots u_{n+1}(1-u_1-\cdots
-u_{n+1}). }
Then, writing $u_0 = 1-u_1-\ldots -u_{n+1}$
\ml{1.9}{\omega_{\Delta_{n+1}} = \frac{du_1\wedge\ldots\wedge du_{n+1}}{u_0\cdot
u_1\cdot u_2\cdots u_{n+1}} = \\
\sum_{i=0}^{n+1} (-1)^{i}d\log(u_0)\wedge\ldots\wedge \widehat{d\log(u_i)}\wedge\ldots\wedge
d\log(u_{n+1}). 
}

Substitute $v_i = tu_i,\ 1\le i\le n+1$. Then $s=1$, and a
calculation   yields
\ga{1.10}{\gamma_{\Delta_{n+1}} =
\frac{1}{v_0}\sum_{i=1}^{n+1}(-1)^i d\log(v_1)\wedge\ldots\wedge
\widehat{d\log(v_i)}\wedge\ldots\wedge d\log(v_{n+1}); \\
\nu_{\Delta_{n+1}} = \sum_{i=0}^{n+1} (-1)^{i}d\log(v_0)\wedge\ldots\wedge \widehat{d\log(v_i)}
\wedge\ldots\wedge d\log(v_{n+1}). \notag
}

For convenience we write
\eq{1.11a}{\nu_{n+1} := \nu_{\Delta_{n+1}}|_{t=0};\quad \gamma_n :=
\gamma_{\Delta_{n+1}}|_{t=0} }

The differential form
$\nu_{n+1}$ 
 plays an important r\^ole in the
computation of de Rham cohomology of the complement of
hyperplane configurations. Let $\sA_t$ be the configuration in
$\P^n$ of $(n+1)$ hyperplanes in general position with
with affine equation $u_0u_1\cdots u_n=0, u_0+u_1+\ldots
+u_n=t\neq 0$. Then $H^n(\P^n\setminus \sA_t)=H^0(\P^n,
\Omega^n_{\P^n}(\log \sA_t))$ is a pure Tate structure 
generated by $\nu_{\Delta{n+1}}$. Let us now make $t$ tend to 0
and consider the degenerate configuration $\sA_0$ with affine
equation $u_0u_1\cdots u_n=0, u_0+u_1+\ldots
+u_n=0$. Exactness of $\nu_{n+1}$ in the case \eqref{1.5}
follows from Aomoto's theory (\cite{A}).

We view $\nu_{\Delta_{n+1}}$ (resp. $\gamma_n$) as a map
$$\sZ_0\Big(\A^{n+1}_k \setminus \{(1-u_1+\cdots+u_{n+1})u_1\cdots
u_{n+1}=0\}\Big)
\to \Omega^{n+1}_k$$ 
(resp. $$\sZ_0\Big(\A^{n+1}_k 
\setminus \{u_1u_2\cdots
u_{n+1}(u_1+\cdots+u_{n+1}) =0\}\Big) \to \Omega^{n}_k).
$$ 
Here $\sZ_0$ denotes the free
abelian group on closed points ($0$-cycles) and the maps are respectively
\eq{1.11}{x \mapsto \text{Tr}_{k(x)/k}\nu|_{\{x\}};\quad x \mapsto
\text{Tr}_{k(x)/k}\gamma|_{\{x\}}. 
}
In the first case, the Nesterenko-Suslin-Totaro theorem
(\cite{S}, \cite{T}) identifies a quotient of
the zero cycles modulo relations coming from curves in $\A^{n+2}$ with the Milnor
$K$-group $K^M_{n+1}(k)$. The evaluation map \eqref{1.11} passes to the quotient, and
the resulting map $K^M_{n+1}(k) \to \Omega_k^{n+1}$ is given on symbols by the
$d\log$-map
\eq{1.12}{ \{x_1,\dotsc,x_{n+1}\} \mapsto d\log(x_1)\wedge\cdots\wedge d\log(x_{n+1}).
}
In the second case, factoring out by curves on $\A^{n+2}$ as in \eqref{1.3a} yields the
Chow group of $0$-cycles $SH^{n+1}(k,n+1)$, and our main result is that evaluation on
$\gamma$ gives an isomorphism
\eq{1.13}{SH^{n+1}(k,n+1) \cong \Omega^n_k.
}

If the hypersurface $X$ in \eqref{1.4a} is smooth, we can also take the residues of
these forms. For a non-zero value of
$t$,
$\text{Res}(\nu_f)$ is a holomorphic $n$-form on $X$. Let $Y^0 \subset Y:F(0,U)=0$ be the
open set of smooth points on the special fibre. We have $\text{Res}(\nu_f|_{t=0}) =
d\text{Res}(\gamma_f|_{t=0})$ on $Y^0$. Evaluation \eqref{1.11} gives a map
\eq{1.14}{CH_0(X) \to \Omega^n_k;\quad x \mapsto
\text{Tr}_{k(x)/k} (\omega|_{x}).
}
What about the special fibre? Consider the special case
\eq{1.15}{f(1,u_1,u_2) = u_1^2 - u_2^3 - au_2 -b;\quad v_1=t^3u_1,\ v_2 = t^2u_2
}
One computes
\eq{1.16}{\nu = \frac{dv_1\wedge dv_2}{v_1^2-v_2^3-at^4v_2-bt^6};\quad \gamma =
\frac{2v_2dv_1 - 3v_1dv_2}{v_1^2-v_2^3-at^4v_2-bt^6}. 
}
One checks that $\text{Res}(\gamma|_{t=0}) = v_2/v_1$ and 
\eq{1.17}{d(v_2/v_1) = -dv_2/2 = d\text{Res}(\nu|_{t=0}) \text{ on }v_1^2-v_2^3=0
}
The assignment $x\mapsto \text{Tr}_{k(x)/k}(v_2/v_1)(x)$ identifies the jacobian of the
special fibre $v_1^2-v_2^3=0$ with $\G_a(k)=k$. 

A final example of specialization, which we understand less well, though it was an
inspiration for this article, concerns the hyperbolic motives of Goncharov \cite{G}. The
{\it matrix coefficients} of his theory (in the sense of \cite{BGSV}) are the objects
$H^{2n-1}(\P^{2n-1}\setminus Q,M \setminus Q\cap M)$. 
Here $Q\subset \P^{2n-1}$ is a smooth quadric and $M$
is a simplex (union of $2n$ hyperplanes in general position).
The subschemes $Q$ and $M$
are taken in general position with respect to each other. The notation is intended to
suggest a sort of abstract relative cohomology group. In de Rham cohomology, the
non-trivial class in $H^{2n-1}_{DR}(\P^{2n-1}\setminus Q)$ is represented by 
\eq{1.18}{\omega = \frac{du_1\wedge\ldots\wedge du_{2n-1}}{(u_1^2+\ldots+u^2_{2n-1}-1)^n}
}
Substituting $u_i=v_it$, we get
\ml{1.19}{\omega = \\
\frac{t^{2n-1}dv_1\wedge\ldots\wedge
dv_{2n-1}+t^{2n-2}dt\wedge\sum (-1)^i
v_idv_1\wedge\ldots\widehat{dv_i}\ldots\wedge
dv_{2n-1}}{\Big(t^2(v_1^2+\ldots+v^2_{2n-1})-1\Big)^n}
}
We get $\nu|_{t=0} = \pm dv_1\wedge\ldots\wedge dv_{2n-1}$ and $\gamma|_{t=0} =
\frac{1}{2n-1}\sum(-1)^iv_idv_1\wedge\ldots\widehat{dv_i}\ldots\wedge
dv_{2n-1}$. Let $\Delta_M \in H_{2n-1}(\P^{2n-1},M; \Z)$ be a generator. The hyperbolic
volume
\eq{1.20}{\int_{\Delta_M}\omega
}
is the {\it real period} (\cite{G} section 4.1) of the Hodge structure associated to
$H^{2n-1}(\P^{2n-1}\setminus Q)$. Goncharov remarks (op. cit.,
Question 6.4 and Theorem 6.5)
that this volume degenerates to the euclidean volume as $t\to 0$. (More precisely, from
\eqref{1.20}, we see that as a relative form, $dv_1\wedge\ldots\wedge dv_{2n-1} = \lim_{t\to 0}
t^{1-2n}\omega$.) He asks for an
interpretation of the degenerated volume in terms of some sort of motive over $k[t]/(t^2)$.

In the Goncharov picture we can view $Q$ as fixed and degenerate $M$. Suppose $M:L_0L_1\cdots
L_{2n-1}=0$ where 
$$L_i = L_i(v_1,\dotsc,v_{2n-1}) = L_i(u_1/t,\dotsc,u_{2n-1}/t)$$
Assuming the $L_i$ are general, clearing denominators and passing to the limit $t\to 0$
yields a degenerate simplex $M_0$ consisting of $2n$ hyperplanes meeting at the point
$v=0$. This limiting configuration leads to the Chow groups $SH^*(k,2n-1)$, but we do
not see how to relate $\gamma|_{t=0}$ to cycles. 

The authors wish to thank A. Goncharov for several very inspiring conversations.

\section{The additive Chow groups}
In this section, we consider  a field $k$, and a $k$-scheme
$X$ of finite type. 

We will throughout use the following notations.
\begin{nota} We set $Q^n={\rm Spec} \ k[t_0,\ldots,
t_n]/(\sum_{i=0}^n t_i)$, together with the faces 
$$\partial_j: Q^{n-1} \to Q^n;\quad \partial_j^*(t_i) = \begin{cases}t_i &
i<j \\ 0 & i=j \\ t_{i-1} & i>j\end{cases}.
$$ 
One also has degeneracies
$$\pi_j:Q^n \to Q^{n-1};\quad \pi_j^*(t_i) = \begin{cases}t_i & i<j \\
t_i+t_{i+1} & i=j \\ t_{i+1} & i>j\end{cases}.
$$
We denote by $\{0\} \in Q^n$ the vertex defined by $t_i=0$. We
write $Q^n_X=Q^n \times_{\Spec(k)} X$. The above face and degeneracy maps make
$Q^\bullet_X$ a cosimplicial scheme. 
\end{nota}

\begin{defn} Let
$\sS\sZ^p(X,n)$ be the free abelian group
on codimension $p$ algebaic cycles on $Q^n_X$ with the property:
\begin{itemize}
\item[(i)] They don't meet $\{0\}\times X$. 
\item[(ii)] They meet all the faces properly, that is in codimension $\ge p$. 
\end{itemize}
\end{defn}
Thus the face maps induce restriction maps
$$\partial_i: \sS\sZ^p(X,n)\to \sS\sZ^p(X, n-1);\quad i=0,\ldots n;\quad \partial =
\sum_{i=0}^n(-1)^i \partial_i;$$  yielding complexes $\sS\sZ^p(X, \bullet)$:
$$
\ldots \stackrel{\partial}{\to} \sS\sZ^p(X, n+1)\stackrel{\partial}{\to} 
\sS\sZ^p(X,n)\stackrel{\partial}{\to} \sS\sZ^p(X,n-1) \stackrel{\partial}{\to} \ldots$$

\begin{defn} The additive Chow groups are given by
$$SH^p(X,n)=H_n (\sS\sZ^p(X, \bullet)).$$ 
\end{defn}

\begin{rmks}\label{rmk2.4}\begin{itemize}\item[(i)] The above should be compared with the
higher Chow groups $CH^p(X,n)$ defined as above with $Q^\bullet$ replaced by
$\Delta^\bullet$, where $\Delta^n := \Spec(k[t_0,\dotsc,t_n]/(\sum t_i -1))$. 
\item[(ii)] The cosimplicial scheme $Q^\bullet$ admits an action of $\G_m$, which we
define by
$$x\star(t_0,\dotsc,t_n) := (t_0/x,\dotsc,t_n/x).
$$
(The reason for the inverse will be clear below.) By functoriality, we obtain a
$\G_m$-action on the $SH^p(X,n)$. 
\item[(iii)] Let $f:X' \to X$ be a proper map with $n=\dim X' -\dim X$. Then
one has a push-forward map $f_*:\sS\sZ^p(X',\bullet) \to \sS\sZ^{p-n}(X,\bullet)$. On
homology this yields $ SH^p(X',n) \to SH^{p-n}(X,n)$. We will be particularly
interested in the case $X=\Spec(k)$, $X' = \Spec(k')$ with $[k':k]<\infty$. We write
$\text{tr}_{k'/k}:SH^p(k',n) \to SH^p(k,n)$ for the resulting map. This trace map
is compatible with the action of $k^\times$ from (ii) in the sense that for $x\in
k^\times$ and $a'\in SH^p(k',n)$ we have $x\star \text{tr}_{k'/k}(a') = 
\text{tr}_{k'/k}(x\star a')$.
\end{itemize}
\end{rmks}

\begin{lem} The $\star$ action of $k$ on $SH^n(k,n)$ is linear, i.e. this action
comes from a $k$-vector space structure on $SH^n(k,n)$. 
\end{lem}
\begin{proof} We extend the action of $k^\times$ on $SH^n(k,n)$ to an action of the
multiplicative monoid $k$ by setting $0\star x=0$.  We have to show that for a closed
point
$x=(u_0,\ldots, u_n) \in Q^n \setminus \cup_{i=0}^n \partial_i(Q^{n-1}) $, and $a,b\in
k$, one has $(a+b)\star x= a\star x + b\star x$. For either $a$ or
$b=0$, this is trivial. Thus we assume $ab\neq 0$. Let $k'=k(x)$. Then the class
in $SH^n(k,n)$ of $x$ is the trace from $k'$ to $k$ of a $k'$-rational point $x'\in
SH^n(k',n)$. Using the compatibility of $\star$ and trace from Remarks
\ref{rmk2.4}(iii) above, we reduce to the case $x$ $k$-rational. 

Let  $\frak{m}=(t_0-u_0, \ldots, t_n-u_n)$ be the sheaf of ideals of $x$. 
Define $\ell(t)= -\frac{ab}{u_0}t+a+b$. Consider the curve $W\subset Q^{n+1}$ defined
parametrically by
$$W = \{(t,-t+\frac{u_0}{\ell(t)}, \frac{u_1}{\ell(t)},\dotsc,\frac{u_n}{\ell(t)})\}.
$$
To check that this parametrized locus is Zariski-closed, we consider the ideal:
$$\sI_W=\Big((t_1+t_0)\ell(t_0)-u_0,t_2\ell(t_0)-u_1,\dotsc,t_n\ell(t_0)-u_n\Big).
$$
If $y=(y_0,\dotsc,y_{n+1})$ is a geometric point in the zero locus of $\sI_W$, then
since the $u_i \neq 0$ we see that $\ell(y_0) \neq 0$. Substituting $t=y_0$, we see that
$y$ lies on the parametrized locus $W$. 

The equation $-t+\frac{u_0}{\ell(t)} = 0$ leads to a quadratic equation in $t$ with
solutions $t=\frac{u_0}{a},\ t=\frac{u_0}{b}$. If $a+b\neq 0$ we have
\ga{}{\partial_0(W) =
(\frac{u_0}{a+b},\dotsc,\frac{u_n}{a+b})=(a+b)\star(u_0,\dotsc,u_n);\notag\\
\partial_1(W) =
(\frac{u_0}{a},\dotsc,\frac{u_n}{a}) + (\frac{u_0}{b},\dotsc,\frac{u_n}{b}) \notag\\ 
= a\star (u_0,\dotsc,u_n) +  b\star (u_0,\dotsc,u_n) \notag \\
\partial_i(W) = 0;\quad i\ge 2, \notag
}
so the lemma follows in this case. If $a+b=0$, then $\partial_0W=0$ as well, and again
the assertion is clear. \end{proof}

\section{Additive Chow groups and Milnor $K$-theory}

We consider the map of complexes $\sZ^p(k, \bullet) \to
\sS\sZ^p(k, \bullet + 1)$ defined on $k$-rational points by 
$\iota: (u_0,\ldots, u_n) \mapsto (-1, u_0,\ldots, u_n)$. This map
induces then a map
\begin{gather}\label{3.1}
\iota: CH^p(k, n)\to SH^{p+1}(k, n+1).
\end{gather} 
By \cite{S} and \cite{T}, one has an isomorphism
\begin{gather}
 K_n^M(k) \cong CH^n(k, n)
\end{gather}
of the higher Chow groups of 0-cycles with Milnor $K$-theory.
It is defined by:
\ga{}{(u_0,\dotsc,u_n) \mapsto \{-\frac{u_0}{u_n},\dotsc,-\frac{u_{n-1}}{u_n}\} \\
\{b_1,\ldots, b_n\} \mapsto (\frac{b_1}{c},\ldots,
\frac{b_n}{c}, -\frac{1}{c});\quad c=-1+\sum_{i=1}^n b_i.
}
Note that if $\sum_{i=1}^n b_i=1$, then the symbol $\beta:=\{b_1,\ldots,
b_n\}$ is trivial in Milnor $K$-theory, and one maps $\beta$ to 0. 

In this way, one obtains a map
\ml{3.5}{K^M_{n-1}(k) \to SH^n(k,n);  \quad
\{x_1,\dotsc,x_{n-1}\} \mapsto \\
\Big(-1,\frac{x_1}{-1+\sum_{i=1}^{n-1}
x_i},\dotsc,\frac{x_{n-1}}{-1+\sum_{i=1}^{n-1} x_i},\frac{-1}{-1+\sum_{i=1}^{n-1}
x_i}\Big). 
}

\section{Differential forms}

In this section we construct a $k$-linear map $\Omega^{n-1}_k \to SH^n(k,n)$. (Here
$\Omega^i_k$ are the absolute K\"ahler differential $i$-forms.) 

The following lemma is closely related to calculations in \cite{C}.
\begin{lem} \label{newdiffform} As a $k$-vector space, 
the differential forms $\Omega^{n-1}_{k}$ are isomorphic to
$(k\otimes _{\Z} \wedge^{n-1} k^\times)/\sR$, where the
$k$-structure on $k\otimes_\Z  \wedge^{n-1} k^\times$ is
via multiplication on the first argument, and where the
relations $\sR$, for $n\geq 2$,  are the $k$-subspace spanned by 
$a\otimes (a\wedge b_1\wedge \ldots \wedge b_{n-2}) +(1-a)
\otimes ((1-a)\wedge b_1 \wedge \ldots \wedge b_{n-2})$, for $
b_i \in k^\times$, $a\in k$,  and where $0\wedge b_2\wedge \wedge
b_{n-2}=0$. 
For $n=1$, one has 
$\sR=\{0\}$. \end{lem}
\begin{proof}
For $n=1$, there is nothing to prove. We assume $n\ge 2$. It will be convenient to
change the relations slightly. Replacing $a$ by $-a$, the relations become
\ml{}{a\otimes (-a\wedge b_1\wedge \ldots \wedge b_{n-2}) -(1+a)
\otimes ((1+a)\wedge b_1 \wedge \ldots \wedge b_{n-2}) \\
=a\otimes (a\wedge b_1\wedge \ldots \wedge b_{n-2}) -(1+a)
\otimes ((1+a)\wedge b_1 \wedge \ldots \wedge b_{n-2}).\notag
}
(To justify this, the $-1$ which appears multiplicatively can be dropped if $k$ has
characteristic $\neq 2$ because the additive group $k$ is $2$-divisible. If $k$ has
characteristic $2$, of course, $-a=a$.)
 One has a
surjective
$k$-linear map
$(k\otimes _{\Z}
\wedge^{n-1} k^\times) \to \Omega^{n-1}_{k/\Z}$ defined by $a\otimes 
(b_1\wedge\ldots \wedge b_{n-1})\mapsto a d\log b_1\wedge \ldots d\log
b_{n-1}$. It factors through $\sR$ as $a d\log a =da=(a+1)d\log (a+1)$
for $a\in k$, with $da=0$ if $a\in \Z$. 

Let us first assume that $n=2$. We define $D: k\to 
(k\otimes _{\Z}  k^\times)/\sR$ by  $D(a)=a\otimes a$ for $a\in
k^\times$, else $D(0)=0$. One has for $a\in k^\times$ the relation
$D(a+b)= (a+b)\otimes (a+b)= a(1+\frac{b}{a})\otimes
a(1+\frac{b}{a})$. By definition of the $k$-structure, this
expression is 
\begin{gather}
=(1+\frac{b}{a})(a\otimes a) +
a\big((1+\frac{b}{a}) \otimes
(1+\frac{b}{a})\big)\notag.\end{gather} 
Modulo
$\sR$, this is 
\begin{gather}
=(1+\frac{b}{a})(a\otimes a) +
a\big((\frac{b}{a})\otimes (\frac{b}{a})\big)=
a\otimes a +\frac{b}{a} (a\otimes a) +
a\big((\frac{b}{a})\otimes (\frac{b}{a})\big). \notag
\end{gather}
Applying again the $k$-structure, one obtains $\frac{b}{a} (a\otimes
a)= b\otimes a$ and 
\begin{gather}
a\big((\frac{b}{a})\otimes (\frac{b}{a})\big)=
b\otimes \frac{b}{a}= b\otimes b - b\otimes a .\notag
\end{gather}
Summing up, one obtains $D(a+b)=D(a)+D(b)$. 
Now if $ab=0$, $D(ab)=0$ and since either $a=0$ or $b=0$ one has
$aD(b)+bD(a)=0.$ Else one
has 
\begin{gather}
-D(ab)+aD(b)+bD(a)= -(ab)\otimes (ab)+ a(b\otimes b)+b(a\otimes a)
=\notag\\
(ab)(1\otimes (ab)^{-1}ab)= (ab)(1\otimes 1)\equiv 0.\notag
\end{gather}
Thus $D$ is a derivation and factors through $\Omega^1_{k} \to (k\otimes 
k^\times)/\sR$ to yield the inverse to the surjection defined above. 
If $n\ge 2$, we extend the map $D$ as follows:
\begin{gather}
D: k^{n-1}\to (k\otimes_\Z \wedge^{n-1}k^\times)/\sR \notag\\
D(b_1,\ldots, b_{n-1}):=D(b_1) \wedge \ldots \wedge D(b_{n-1}):=
b_1\cdot \ldots \cdot b_{n-1} \otimes b_1 \wedge \ldots
b_{n-1}.\notag \end{gather} 
This symbol is immediately seen to be alternating. 
Furthermore, the computations above yields 
\begin{gather}
D(ab)\wedge D(b_2)\wedge \ldots \wedge D(b_{n-1}) =\notag\\
aD(b)\wedge D(b_2)
\wedge \ldots \wedge D(b_{n-1}) + bD(a)\wedge D(b_2)
\wedge \ldots \wedge D(b_{n-1}),\notag
\end{gather}
and
\begin{gather}
D(a+b)\wedge D(b_2) \ldots \wedge D(b_{n-1})=\notag \\ 
D(a)\wedge D(b_2) \ldots \wedge D(b_{n-1})+
D(b)\wedge D(b_2) \ldots \wedge D(b_{n-1}).\notag
\end{gather}
This defines the inversed map to the surjection 
$(k\otimes_\Z \wedge^{n-1}k^\times)/\sR 
\to \Omega^{n-1}_{k}$ defined above and proves the lemma. 
\end{proof}

\begin{prop} \label{inv}
One has a well-defined $k$-linear map
\begin{gather}
\phi: \Omega^{n-1}_{k} \to SH^n(k,n)\notag\\
\alpha:=a d\log b_1\wedge \ldots \wedge d\log b_{n-1} \mapsto a\star (-1, 
\frac{b_1}{c},\ldots, \frac{b_{n-1}}{c}, -\frac{1}{c})\notag
\end{gather}
where $c=-1+\sum_{i=1}^{n-1}b_i$. The diagram
$$\begin{CD} K^M_{n-1}(k) @>d\log >> \Omega^{n-1}_k \\
@V\iota V \eqref{3.1}V @V\phi VV \\
SH^n(k,n) @= SH^n(k,n)
\end{CD}
$$
is commutative. 
\end{prop}
\begin{proof} If $c=-1+\sum_{i=1}^{n-1} b_i=0$, then $\alpha =0$,
else $a\star (-1, 
\frac{b_1}{c},\ldots, \frac{b_{n-1}}{c}, -\frac{1}{c})$ is
defined. We now change notation and write $c=-1+a +\sum_{i=2}^{n-1} b_i$.  By Lemma
\ref{newdiffform}, we have to show
\ml{}{
0=\rho:= \\
a\star (-1, 
\frac{a}{c},\ldots, \frac{b_{n-1}}{c}, -\frac{1}{c}) 
- (a+1)\star (-1, 
\frac{a+1}{c+1},\ldots, \frac{b_{n-1}}{c+1}, -\frac{1}{c+1}). \notag
}

If $a=0$, then one has
\begin{gather}\rho=- (-1, \frac{1}{c},\ldots,
\frac{b_{n-1}}{c}, -\frac{1}{c}) =
- \iota\{1,b_2,\ldots, b_{n-1}\}=0.\notag
\end{gather}
Similarly, $\rho=0$ if $a=-1$. Assume now $a\neq 0,-1$. 
Set $b=-a \in k\setminus \{0,1\}$. One defines, for $n\ge 3$ and
$(u_1,\ldots, u_{n-1}) \in (\Delta^{n-2}\setminus
\cup_{i=1}^{n-1} \Delta^{n-3})(k)$, the parametrized curve
\begin{gather}
\Gamma(b,u):=\{(\frac{-1}{b}+t, \frac{1}{b-1}, -t,
\frac{-u_1}{b(b-1)},\ldots, \frac{-u_{n-1}}{b(b-1)})\} \subset
Q^{n+1} \notag  \end{gather}
and for $n=2$
\begin{gather}  
\Gamma(b):=\{(\frac{-1}{b}+t, \frac{1}{b-1}, -t,
\frac{-1}{b(b-1)})\}\subset Q^3.
\end{gather}
(See \cite{T} for the origin of this definition). This curve
is indeed in good position, so it lies in $\sS\sZ^n (Q^{n+1})$.
One computes 
\begin{gather}\label{4.2}
\partial \Gamma(b, u)=(1-b)\star 
(-1, 1-\frac{1}{b}, \frac{u_1}{b},\ldots,
\frac{u_{n-1}}{b})\\ + 
b\star (-1, \frac{b}{b-1}, \frac{-u_1}{b-1},\ldots,
\frac{-u_{n-1}}{b-1}).\notag
\end{gather}
(Resp. in the case $n=2$
$$\partial\Gamma(b) = (1-b)\star (-1,1-\frac{1}{b},\frac{1}{b})
+b\star(-1,\frac{b}{b-1},\frac{-1}{b-1}).)
$$
Now one has 
\begin{gather}
(-1, 1-\frac{1}{b}, \frac{u_1}{b},\ldots,
\frac{u_{n-1}}{b})=\iota 
\{\frac{1-b}{u_{n-1}}, -\frac{u_1}{u_{n-1}},\ldots,
-\frac{u_{n-2}}{u_{n-1}}\} \notag\\
= \iota[\{1-b,-\frac{u_1}{u_{n-1}},\ldots,
-\frac{u_{n-2}}{u_{n-1}}\} -\{u_{n-1},u_1,\ldots,
u_{n-2}\}],\notag 
\end{gather}
as the rest of the multilinear expansion contains only symbols
of the shape $\{\ldots, u_{n-1},\ldots, -u_{n-1}, \ldots\}$. 
On the other hand, since
$\sum_{i=1}^{n-1}u_i=1$,
one has 
$\{u_{n-1},u_1,\ldots,
u_{n-2}\}=0$.  Similarly, one has
\begin{gather}
(-1, \frac{b}{b-1}, \frac{-u_1}{b-1},\ldots,
\frac{-u_{n-1}}{b-1})=\iota
\{\frac{b}{u_{n-1}}, -\frac{u_1}{u_{n-1}},\ldots,
-\frac{u_{n-2}}{u_{n-1}}\}.\notag
\end{gather}
The same argument yields that this is
\begin{gather} 
\iota\{b,-\frac{u_1}{u_{n-1}},\ldots,
-\frac{u_{n-2}}{u_{n-1}}\}.\notag
\end{gather}

It follows now from \eqref{4.2} that for $n\ge 3$ we have the relation in $SH^n(k,n)$
\ml{}{(1-b)\star \iota\{1-b,-\frac{u_1}{u_{n-1}},\ldots,
-\frac{u_{n-2}}{u_{n-1}}\} + \\
b\star\iota\{b,-\frac{u_1}{u_{n-1}},\ldots,
-\frac{u_{n-2}}{u_{n-1}}\} = 0
}
(The analogous relation 
for $n=2$ is similar.) The assertion of the proposition now follows
from Lemma \ref{newdiffform}. 
\end{proof}
\begin{rmk} \label{coord}
The coordinates $u_i$ of $Q^n$ respecting
the boundary $\partial(Q^{n+1})$ are defined up to scalar
$\lambda \in k^\times$, thus $\phi_{\lambda u_i}(\alpha)=
\lambda \star \phi_{u_i}(\alpha)$, and the map $\phi$ is
coordinate dependent. 
\end{rmk}
\begin{prop}With notation as above, the map 
$$\phi: \Omega^{n-1}_k \to SH^n(k,n)$$ 
is surjective. In particular, $SH^n(k,n)$ is generated by the classes of $k$-rational
points in $Q^n$.
\end{prop}
\begin{proof} It is easy to check that the image of $\phi$ coincides with the subgroup
of $SH^n(k,n)$ generated by $k$-points. Clearly, $SH^n(k,n)$ is generated by closed
points, and any closed point is the trace of a $k'$-rational point for some finite
extension $k'/k$. We first reduce to the case $k'/k$ separable. If $x\in Q^n_k$ is a
closed point in good position (i.e. not lying on any face) such that $k(x)/k$ is not
separable, then a simple Bertini argument shows there exists a curve $C$ in good
position on $Q^{n+1}$ such that
$\partial C = x + y$ where $y$ is a zero cycle supported on points with separable
residue field extensions over $k$. Indeed, let $W\subset Q^{n+1}$ be the union of the
faces. View $x\in W$. Since $x$ is in good position, it is a smooth point of $W$.
Bertini will say that a non-empty open set in the parameter space of $n$-fold
intersections of hypersurfaces of large degree containing $x$ will meet $W$ in $x$ plus
a smooth residual scheme. Since $k$ is necessarily infinite, there will be such an
$n$-fold intersection defined over $k$. Since the residual scheme is smooth, it cannot
contain inseparable points. Then $x\equiv -y$ which is supported on separable points. 

We assume now $k'/k$ finite separable, and we must show that the trace of a $k'$-point
is equivalent to a zero cycle supported on $k$-points.  
Since the image of $\phi$ is precisely the subgroup generated by $k$-rational points, it
suffices to check that the diagram
\eq{4.4}{\begin{CD} \Omega^{n-1}_{k'} @>\phi>> SH^n(k',n) \\
@VV\text{Tr}_{k'/k}V @VV\text{Tr}_{k'/k}V \\
\Omega^{n-1}_{k} @>\phi>> SH^n(k,n) 
\end{CD}
}
commutes. Because $k'/k$ is separable, one has $\Omega^{n-1}_k \inj \Omega^{n-1}_{k'}$,
and $\Omega^{n-1}_{k'} = k'\cdot \Omega^{n-1}_{k}$. One reduces to showing, for $\alpha =
(-1,\alpha_1,\dotsc,\alpha_n)\in Q^n(k)$ and $t\in k'$, that $\text{Tr}(t\star \alpha) =
(\text{Tr}(t))\star \alpha$.

Let $P(V)=V^N+a_{N-1}V^{N-1}+\ldots +a_1V+a_0 \in k[V]$ be the minimal
polynomial of $-\frac{1}{t}$. We set $b_N=\frac{-1}{\alpha_n}$, 
$b_i=\frac{-a_i}{\alpha_n}, i=N-1,\ldots, 2$ and $b_i=a_i, i=1,0$.
We define the polynomial
$Q(V,u)=b_NV^{N-1} u + \ldots b_2Vu + b_1V + b_0 \in k[V,u]$,
which by definition fulfills $Q(V,-\alpha_nV)=P(V)$. 
We define the ideal 
\begin{gather}
\sI =(Q(V_0,u), V_1+\alpha_1 V_0, \ldots, V_{n-1} +\alpha_{n-1}
V_0)\notag \\
\subset k[V_0,\ldots,V_{n-1}, u].\notag
\end{gather}
It defines a curve $W\subset \A^{n+1}$. We think of $\A^{n+1}$
as being $Q^{n+1}$ with the faces $V_0=0,\ldots, V_{n-1}=0, u=0,
u+\sum_{i=0}^{n-1} V_i =0$.  Then this curve is in general
position and defines a cycle in $\sS\sZ^1(k, n+1)$.

Since $b_0\neq 0$, and $\alpha_i \neq 0$, one has 
\begin{gather}
\partial_i W=0 , i=0,1,\ldots, n-1.
\end{gather}
One has
\begin{gather}
\partial_uW \ \text{defined \ by} \ (a_1V_0 +a_0V_0,
V_1+\alpha_1 V_0,\dotsc,V_{n-1}+\alpha_{n-1}V_0).\notag
\end{gather}
To compute the last face, we observe that the ideal
$$(u+
\sum_{i=0}^{n-1} V_i =0,V_1+\alpha_1 V_0, \ldots, V_{n-1} +\alpha_{n-1}
V_0), $$
contains $u+\alpha_n V_0$. Consequently $\partial_{u+
\sum_{i=0}^{n-1} V_i}W$ is defined by 
\begin{gather}
\Big(Q(V_0,-\alpha_nV_0)=P(V_0), V_1+\alpha_1 V_0, \ldots, V_{n-1} +\alpha_{n-1}
V_0\Big).\notag \end{gather}
Thus one obtains
\begin{gather}
0\equiv (-1)^{n}\partial W= \frac{a_1}{a_0} \star(-1,\alpha)
- t\star(-1, \alpha). \notag\end{gather}
Since $P$ is the minimal polynomial of $-\frac{1}{t}$, 
$\frac{a_1}{a_0}$  is the trace of $t$. 
\end{proof}

\section{The main theorem}

Recall \eqref{1.11a} we have a logarithmic $(n-1)$ form $\gamma_{n-1}$ on \break $Q^n =
\Spec(k[v_0,\dotsc,v_n]/(\sum v_i))$
\ga{5.1}{\gamma_{n-1} = \frac{1}{v_0}\sum_1^n (-1)^i d\log(v_1)\wedge\ldots\wedge
\widehat{d\log(v_i)}\wedge\ldots\wedge d\log(v_n) \\
d\gamma_{n-1} = \nu_n = \sum_0^n (-1)^i d\log(v_0)\wedge\ldots\wedge
\widehat{d\log(v_i)}\wedge\ldots\wedge d\log(v_n) \notag
}
Writing $v_i = V_i/V_{n+1}$, we can view $\gamma_{n-1}$ as a meromorphic form on $\P^n =
\text{Proj}(k[V_0,\dotsc,V_{n+1}]/(\sum_0^n V_i))$. Let $\sA:V_0\cdots V_n=0;\
\infty:V_{n+1}=0$. The fact that
$d\gamma_{n-1}$ has log poles on the divisors $V_i=0,\ 0\le i\le
n$ implies that  
\eq{5.2}{\gamma_{n-1} \in \Gamma\Big(\P^n,\Omega^{n-1}_{\P^n}
(\log(\sA+\infty))(-\infty)\Big)
}
In particular, $\gamma_{n-1}$ has log poles, so we can take the
residue  along components of $\sA$. 
(The configuration $\sA$ does not have normal crossings.
The sheaf $\Omega^{n-1}_{\P^n}
(\log(\sA+\infty))$ is defined to be
the subsheaf of $j_*\Omega^{n-1}_{\P^n-\sA-\infty}$ generated
by forms without poles and the evident log forms with residue $1$
along one hyperplane and $(-1)$ along another one. According to
\cite{A}, the global
sections of this naturally defined log sheaf compute
de Rham cohomology. )

\begin{lem}\label{lem5.1} We have the following residue formulae
\ga{5.3}{{\rm Res}_{v_i=0} \gamma_n = (-1)^i \gamma_{n-1};\quad 1\le i\le n+1\\
{\rm Res}_{v_0=0} \gamma_n = \gamma_{n-1}.\notag
}
\end{lem}
\begin{proof}The assertion for $1\le i\le n+1$ is straightforward. For $i=0$, one can
either compute directly or argue indirectly as follows:
\eq{5.4}{d\text{Res}_{v_0=0}\gamma_n = \text{Res}_{v_0=0}d\gamma_n =
\text{Res}_{v_0=0}\nu_{n+1} = \nu_n = d\gamma_{n-1}. 
}
Since the global sections \eqref{5.2} has dimension $1$, this suffices to determine
$\text{Res}_{v_0=0}\gamma_n$. (To verify dim. $1$, let $\sA' \subset \sA$ be defined by
$V_1\cdots V_n=0$. Then $\sA'+\infty$ consists of $n+1$ hyperplanes in general position
in $\P^n$, so $\Omega^{n-1}_{\P^n}(\log(\sA'+\infty)) =
\wedge^{n-1}\Omega^{1}_{\P^n}(\log(\sA'+\infty))\cong
\sO^{\oplus n}_{\P^n}$. One looks at the evident
 residue 
$$\Omega^{n-1}_{\P^n}
(\log(\sA+\infty))(-\infty) \to \Omega^{n-2}_{\P^{n-1}}
(\log(\sA+\infty))(-\infty).
$$
along $V_0=0$.) 
\end{proof}
\begin{thm} \label{mainthm} The assignment $x \mapsto {\rm Tr}_{k(x)/k}(\gamma(x))$ gives an isomorphism
$$SH^n(k,n) \cong \Omega^{n-1}_k.
$$
\end{thm}
\begin{proof} Let $X\subset Q^{n+1}_k$ be a curve in good position. For a zero-cycle $c$
in good position on $Q^n$ we write $\text{Tr}\gamma_{n-1}(c)\in \Omega^{n-1}_k$
(absolute differentials) for the evident linear combination of traces from residue fields
of closed points. We must show $\text{Tr}\gamma_{n-1}(\partial X)=0$. Let $\overline{X}$
denote the closure of $X$ in $\P^n$. The form $\gamma_{n-1}$ dies when restricted to
$\infty$. We may therefore replace $\partial X$ with $\partial \overline{X}$.

Consider the diagram with $D_j:= \partial_j(Q^n)\cap
\overline{X},\ D= \sum_{j=0}^n D_j$
\eq{5.5}{\begin{CD}H^0(\overline{X},\Omega^n_{\overline{X}/\Z}(\log D)) @>\sum
\text{Res}_{D_j}>>
\oplus_{j=0}^n
\Omega^{n-1}_{D_j} \\
@VVV @VV a V \\
H^0(\overline{X}\setminus D,\Omega^n_{\overline{X}/\Z}) @>\delta >> \oplus_{j=0}^n
H^1_{D_j}(\overline{X},\Omega^n_{\overline{X}/\Z}) @>b>>
H^1(\overline{X},\Omega^n_{\overline{X}/\Z}) \\
@. @. @VV\text{Tr} V \\
@. @. \Omega^{n-1}_k.
\end{CD}
}

The map $\text{Tr}\circ b\circ a$ is the trace map used to
define our map. By Lemma 
\ref{lem5.1}, $\gamma_{n-1}(D_j) = (\text{Res}_{v_j=0}\gamma_n)(D_j)$. The desired
vanishing follows, taking $\gamma_n \in H^0(\overline{X},\Omega^n_{\overline{X}/\Z}(\log
D))$ and using $b\circ \delta = 0$.

We now have
$$\Omega^{n-1}_k \stackrel{\phi}{\surj} SH^n(k,n)
\stackrel{\gamma_{n-1}}{\longrightarrow}
\Omega^{n-1}_k.
$$
It suffices to check the composition is multiplication by $(-1)^{n+1}$. Given
$b_1,\dotsc,b_{n-1}
\in k$ with $c:=\sum b_i -1\neq 0$, the composition is computed
to be (use \eqref{5.1}
and Proposition \ref{inv})
\ml{}{a\cdot d\log(b_1)\wedge\ldots\wedge d\log(b_{n-1}) \mapsto \\
 -a\sum (-1)^i
d\log(\frac{b_1}{ac})\wedge\ldots\wedge
\widehat{d\log(\frac{b_i}{ac})}\wedge\ldots\wedge d\log(\frac{-1}{ac}) \notag
}
Expanding the term on the right yields
\ml{}{ -a\sum (-1)^i
d\log(\frac{b_1}{c})\wedge\ldots\wedge
\widehat{d\log(\frac{b_i}{c})}\wedge\ldots\wedge d\log(\frac{-1}{c}) + \\
-a\cdot d\log(a)\wedge\Big(\ldots\Big), \notag
}
and it is easy to check that the terms involving $d\log(a)$ cancel. In this way, one
reduces to the case $a=1$. Here
\ml{}{-\sum (-1)^i
d\log(\frac{b_1}{c})\wedge\ldots\wedge
\widehat{d\log(\frac{b_i}{c})}\wedge\ldots\wedge d\log(\frac{-1}{c}) = \\
(-1)^{n+1} \frac{db_1}{b_1}\wedge\ldots
\wedge \frac{db_{n-1}}{b_{n-1}} + \frac{dc}{c}\wedge\Big(\ldots\Big).
}
Again the terms involving $d\log(c)$ cancel formally, completing the proof. 
\end{proof}
\begin{rmk} The isomorphism of the main theorem \ref{mainthm} depends,
according to the Remark \ref{coord}, on the scale of
the coordinates. 
\end{rmk}

\begin{chal}
Finally, as a challenge we remark that the K\"ahler differentials have operations 
(exterior derivative, wedge product,...) which are not evident on the cycles $\sS\sZ$.
For  example, one can show that the map
\begin{gather}
\nabla(x_0,\ldots, x_n)= (x_0 , -\frac{x_1x_0}{1-x_0},\ldots,-\frac{x_nx_0}{1-x_0},
-\frac{x_0}{1-x_0}) \end{gather}
satisfies
\begin{gather}
\gamma_n (\nabla(x))=(-1)^nd\gamma_{n-1}(x)
\end{gather}
and hence induces the exterior derivative on the $0$-cycles. The map is not uniquely
determined by this property, and this particular map does not preserve good
position for cycles of dimension $>0$. Can one find a  geometric correspondence on the
complex
$\sS\sZ^\bullet$ which induces $d$ on the $0$-cycles? What about the pairings $(a,b)
\mapsto a\wedge b$ or $(a,b) \mapsto a\wedge db$? 
\end{chal}

\newpage
\bibliographystyle{plain}

\begin{thebibliography}{99}
\bibitem{V} Voevodsky, V.: Motivic cohomology groups are
isomorphic to higher Chow groups in any characteristic,
$K$-theory server.
\end{thebibliography}
\renewcommand\refname{References}

\end{document}